\documentclass[twoside,12pt,a4paper]{article}
\usepackage[marked,printedin,
]{myart2k}

\ppnum{}{\htmladdnormallink{funct-an/9704001}%
  {http://arXiv.org/abs/funct-an/9704001}}{1997}
\def\ifundefined#1{\expandafter\ifx\csname#1\endcsname\relax}
\providecommand{\comment}[1]{}

\usepackage{amsfonts}
\providecommand{\tthdump}[1]{#1}
\tthdump{}
\newcommand{\algebra}[1]{\ensuremath{\mathfrak{#1}}}

\newcommand{\Cliff}[2][\comment]{{\ensuremath{%
\mathcal{C}\kern-0.18em\ell(#1,#2)}}}
\newcommand{\object}[2][\,]{\ensuremath{\mathrm{#2}#1}}
\newcommand{\Space}[3][\relax]{\ensuremath{\mathbb{#2}^{#3}_{#1}}}

\newcommand{\FSpace}[3][\relax]{\ensuremath{ #2_{#3}^{#1}}}

\newcommand{\bos}{\mathcal{B}}
\ifundefined{qed}
    \DeclareMathSymbol{\qed}{0}{AMSa}{"03}
\fi
\newcommand{\modulus}[2][\relax]{\left| #2 \right|_{#1}}
\newcommand{\scalar}[3][\relax]{\left\langle #2,#3 \right\rangle_{#1}}
\newcommand{\intersect}[2]{#1\mid_{#2}}
\newcommand{\relstack}[2]{{\mathrel {\mathop {#1}\limits _{#2}}}}
\providecommand{\eqref}[1]{\textup{(\ref{#1})}}

\newcommand{\person}[1]{\textsc{#1}}
\providecommand{\href}[2]{#2}

\input{cyr.fd}
\makeatletter
\def\item{\@inmatherr \item \@ifnextchar [\@item %
  {\@noitemargtrue \@item [\textup{\@itemlabel} ]}}
\def\p@enumi{\thethm.}
\makeatother
\renewcommand{\FSpace}[2]{\ensuremath{\mathcal{#1}_{#2}}}
\newcommand{\bin}[2]{{#1 \choose #2}}
\newcommand{\oper}[1]{\mathcal{#1}}

\begin{document}

\title[Umbral Calculus and Cancellative Semigroups]{Umbral Calculus\\
  and Cancellative Semigroup Algebras\thanks{Supported
    by grant 3GP03196 of the FWO-Vlaanderen (Fund of Scientific 
    Research-Flanders),
    Scientific Research Network ``Fundamental Methods and Technique in
    Mathematics''.}}

\author[Vladimir V. Kisil]
{\href{http://amsta.leeds.ac.uk/kisilv~/}{Vladimir V. Kisil}}
\address{Department of Pure Mathematics\\
  University of Leeds\\
  Leeds LS2\,9JT\\
  UK
  }
\email{\href{mailto:kisilv@amsta.leeds.ac.uk}%
  {kisilv@amsta.leeds.ac.uk}}
\urladdr{\href{http://amsta.leeds.ac.uk/~kisilv/}%
  {http://amsta.leeds.ac.uk/\~{}kisilv/}}
\date{April 1, 1997; Revised February 4, 2000}
\maketitle
\dedicatory{Dedicated to the memory of Gian-Carlo Rota}
\begin{abstract}
  We describe some connections between three different fields:
  combinatorics (umbral calculus), functional analysis (linear
  functionals and operators) and harmonic analysis (convolutions on
  group-like structures). Systematic usage of cancellative semigroup,
  their convolution algebras, and \emph{tokens} between them provides a
  common language for description of objects from these three fields.
  \AMSMSC{43A20}{05A40, 20N++} 
\end{abstract}
\tableofcontents
\bigskip

\epigraph{It is a common occurrence in combinatorics for a simple idea
  to have vast ramifications; the trick is to realize that many
  seemingly disparate phenomena fit under the same umbrella.}%
{Richard P. Stanley,~\cite{Stanley95a}.}{}

\section{Introduction}
We discuss some interactions between representation theory of
semigroups, umbral calculus in combinatorics, and linear
operators in functional analysis.  Our main bridge between them is the
notion of cancellative semigroups, convolution algebras over them and
tokens from such a semigroup to another.

We look for a vocabulary, which is reasonably general
and allows us to translate without alterations as much about umbral
calculus as possible. The language we select relies on convolution
algebras. Such algebras arise from essentially different sources like
groups, posets, or linear operators. The goal of the present paper to
study properties of convolutions algebras, which are independent from
their origins.
 
The \emph{umbral calculus} was known long ago but it was not well
understood till very recent time. During its history it took many
different faces: linear functionals~\cite{Rota64a}, Hopf
algebras~~\cite{RomRota78}, axiomatic definition~\cite{RotaTay94a} to
list only few different interpretations. 
\comment{
  It was born as a sort of mysterious
  \emph{witchcraft}. To became a rigorous theory it firstly take a
  clothes of \emph{functional analysis}~\cite{RotaTay94a}:
  
  \begin{quotation}
    {\ldots} it was one of the present authors
    [\person{Rota}~\cite{Rota64a}] who first pointed out explicitly what
    nowadays had become completely obvious, namely, that the shift from
    $a^n$ to $a_n$ can be ``explained'' as the application of a linear
    functional on the space of polynomials in the variable $a$.
  \end{quotation}
  
  However it was only beginning of complicated process and more involved
  \emph{Hopf algebras} was introduced~\cite{RotaTay94a}:
  
  \begin{quotation}
    \ldots it took the same author another twenty years to realize that
    this simple truth could not alone explain the stunts of calculations
    performed by the umbral mathematics. At long last, it was realized
    that umbral calculus could be made entirely rigorous by using the
    language of Hopf algebras, and this was done in a lengthy
    treatment~\cite{RomRota78}.
  \end{quotation}
  
  And this was not the final step~\cite{RotaTay94a} yet:
  
  \begin{quotation}
    However, although the notation of Hopf algebra satisfied the most
    ardent advocate of spic-and-span rigor, the translation of
    ``classic'' umbral calculus into the newly found rigorous language
    made the method altogether unwieldy and unmanageable.
  \end{quotation}
  
  Thus in the recent paper~\cite{RotaTay94a} a new formulation of the
  umbral calculus was described. The main feature of the
  semantic-axiomatic approach is an explicit separation of three
  equivalence relations (denotes by $=$, $\simeq$, $\equiv$) and two
  multiplication by integers ($n\cdot \alpha$ and $n.\alpha$).
}

The goal of our paper is to give a model, which realizes described
features by means of common objects. What is a porpoise of such
models? One can think on the Poincar\'e model of non-Euclidean
geometry in a disk. The model does not only demonstrate the logical
consistency of Lobachevski's geometry. Another (no less important)
psychological function is to express the unusual geometry via objects
of the intuitive classic geometry. The logical consistence of umbral
calculus is evident (from papers~\cite{RomRota78,Rota64a,RotaTay94a}
for example). But is our feeling that some more links with fundamental
objects are still desirable.

The paper presents basic definitions and properties of cancellative
semigroups, tokens, convolution algebras, and their relations with
umbral calculus. We will illustrate our consideration mostly by the
simplest example: the cancellative semigroup of non-negative integers
$\Space[+]{N}{}$. More non-trivial examples will be present elsewhere.

The layout is as follows. In the next Section we sketch three
different interpretations of umbral calculus from
papers~\cite{RomRota78,Rota64a,RotaTay94a}.
Section~\ref{sec:canc-semigr-tokens} introduces cancellative
semigroups, convolutions, and tokens and their basic
properties. We use them in Section~\ref{sec:shift-invar-oper} to
describe principal combinatorial objects like delta families, generating
functions, and recurrence operators. We conclude the paper with
Section~\ref{sec:models-umbr-calc} which links our construction with
three realizations of umbral calculus
from~\cite{RomRota78,Rota64a,RotaTay94a} recalled in
Section~\ref{se:umbral}.

\section{The Umbral Calculus}\label{se:umbral}

The umbral calculus was put on the solid ground by \person{G.-C.~Rota}
and collaborates~\cite{RomRota78,Rota64a,Rota95,RotaTay94a} in three
different ways at least. We repeat here some essential definitions
from these papers.

\subsection{Finite Operator Description}\label{ss:operator-d}
The three ways is based essentially on two fundamental
notions~\cite{Rota78}:
\begin{enumerate}
\item A \emph{polynomial sequence of binomial type}, that is, a
  sequence of polynomials $p_n(x)$ ($\deg p_n =n$) with complex
  coefficients, satisfying the identities:
  \begin{equation}\label{eq:binom}
    p_n(x+y)=\sum_{k=0}^{n} p_k(x) p_{n-k}(y).
  \end{equation} 
  Our definition is different from the original
  one---we take polynomial divided by $n!$. Our choice will be
  explained shortly in connection with the group property. Such form
  of the defining identity together with a more accurate name
  \emph{polynomial sequence of integral type} were used
  in~\cite{BarnabBriniNicol80}. Paper~\cite{Knuth99a} also start from
  this formula, the polynomials are called there \emph{convolution
  polynomials} and their investigation is made with a help of
  \textsl{Mathematica} software.

\item A \emph{shift invariant operator}, namely a linear operator $S$
  on the vector space $\FSpace{P}{}$ of all such polynomials, which
  commutes with the ordinary derivative $Dp(x)=p'(x)$, that is, an
  operator $S$ with the property that $TDp(x)=DTp(x)$ for all
  polynomials $p(x)$. Within shift invariant operators the following
  one plays an exceptional role: an operator $Q$ is said to be a
  \emph{delta operator} (\emph{associated} to a polynomial sequence $p_n(x)$)
  when
  \begin{equation}\label{eq:l-diff}
    Q p_n(x)= p_{n-1} (x) \textrm{ or more general } Q^k p_n(x)=
    p_{n-k}(x).
  \end{equation}
\end{enumerate}

\subsection{Hopf Algebras Description}
There is a canonical isomorphism between shift invariant operators
and linear functionals (see Proposition~\ref{pr:isomorphism}), thus we
can think on powers of delta $Q^k$ as functionals $q_k$ on
polynomials connected with associated polynomial sequence $p_n$ by the
\emph{duality}:
\begin{equation}
  \scalar{q_k}{p_n(x)}=\delta _{kn},
\end{equation} 
which easily follows from~\eqref{eq:l-diff}. Here
$\delta_{kn}$ is the Kronecker delta\footnote{We use the same letter
  $\delta$ to denote also the Dirac delta function. The Kronecker
  delta $\delta_{kn}$ and the Dirac function $\delta(x)$ can be
  distinguished by their arguments.}.  We are going to show that
delta operators and polynomial sequences is not only the dual notions,
but two faces of the same object---tokens between
cancellative semigroups (Examples~\ref{ex:polin1} and~\ref{ex:oper}).

\subsection{The Semantic Description}
This subsection contains definitions
from~\cite{RotaTay93a,RotaTay94a}.

\begin{defn}\cite{RotaTay93a}
  An \emph{umbral calculus} is an ordered quadruple, consisting of
\begin{enumerate}
  
\item A commutative integral domain $D$ whose quotient field is of
  characteristic zero.
  
\item An alphabet $A$ whose elements are termed \emph{umbrae} and
  denoted by Greek letters. This alphabet generates the polynomial
  ring $D[A]$, elements of which are termed \emph{umbral polynomials}.
  
\item A linear functional $\object{eval}: D[A] \rightarrow D$, such
  that $\object{eval}(1)=1$, i.e., $\object{eval}$ leaves $D$ fixed,
  and such that
\begin{equation}\label{eq:eval-sep}
  \object{eval}(\alpha^i \beta^j \cdots \gamma^k)=
  \object{eval}(\alpha^i) \object{eval}(\beta^j) \cdots \object{eval}
  (\gamma^k)
\end{equation}
where $\alpha $, $\beta$, \dots, $\gamma$ are distinct umbrae

\item\label{it:augmentation} A distinguished umbrae $\epsilon$ in $A$,
  such that $\object{eval}(\epsilon^i)=\delta _{i,0}$. $\epsilon$ is
  sometimes called the augmentation.
\end{enumerate}
\end{defn}

\begin{defn}\cite{RotaTay93a}
  Let $f$ and $g$ be two umbral polynomials in $D[A]$.  Then $f$ and
  $g$ are said to be \emph{umbrally equivalent}, written $f\simeq g$
  when $\object{eval} (f) = \object{eval}(g)$.
\end{defn}

\begin{defn}\cite{RotaTay93a}
  Two umbral polynomials $p$ and $q$ are termed \emph{exchangeable},
  written $p\equiv q$, when $p^n \simeq q^n$ for all $n\geq 0$.
\end{defn}

\section{Cancellative Semigroups and Tokens}\label{sec:canc-semigr-tokens}
\subsection{Cancellative Semigroups: Definition and Examples}
Our consideration is based on the following
notion of \emph{semigroups with cancellation} or \emph{cancellative
  semigroups}~\cite[\S~IV.1.1]{AlgebraII}: 
\begin{defn}
  A semigroup $C$ is called a \emph{left (right) cancellative
    semigroup} if for any $a$, $b$, $c\in C$ the identity $ca
  = cb$ ($ac=bc$) implies $a=b$. A \emph{cancellative semigroup} is
  both left and right cancellative semigroup.
\end{defn} 
Equivalently a left (right) cancellative semigroup is defined
by the condition that for arbitrary $a,b\in C$ the equation $a\cdot
x=b$ ($x\cdot a=b$) has at most one solution (if any). We hope that
readability of the paper will be better if we use ``\emph{c-semigroup}''
to denote ``cancellative semigroup''. We also use the notion of
\emph{c-set} $C$ (which is weaker than c-semigroup) by allowing
the multiplication $(a,b) \mapsto ab$ on $C$ be defined only on a
proper subset of $C \times C$.
\begin{defn}
  An element $e$ is called a (left,right) \emph{source} of a (left,
  right) c-semigroup $C$ if both equations (the first equation, the
  second equations) $x\cdot e=b$, $e\cdot x=b$ do (does, does) have
  a (unique) solution for any $b\in C$.
\end{defn}
Obviously, a (left,right) identity on $C$ will be a (left,right)
source. 

We still denote the unique-if-exist solutions to equations $x \cdot
a=b$ and $a \cdot x=b$ by $[ba^{-1}]$ and $[a^{-1}b]$
correspondingly. Here the braces stress that both $[ba^{-1}]$ and
$[a^{-1}b]$ are \emph{mono}symbols and just ``$a^{-1}$'' is not
defined in general.  Let us also assume that all c-semigroups under
consideration can be equipped by an \emph{invariant measure} $db$,
namely $db=d(ab)$ for all $a\in C$. We will not discuss herein
conditions for its existence or possible modifications of our
constructions for the case of \emph{quasi}-invariant measures.


C-semigroups were investigated as an algebraic object in connection
with groups. It is particularly known that any commutative c-semigroup
can be embedded in a group but this is not necessary true for a
non-commutative c-semigroup~\cite[\S~IV.1.1]{AlgebraII}. Our
motivation in this object is the following.  If we have a group $G$
then the important associated object is the left regular (linear)
representation by ``shifts'' $\pi_g f(h)=f(hg)$ in
$\FSpace{L}{p}(G,d\mu)$. Then one can introduce their linear
span---convolutions:
\begin{equation}\label{eq:conv}
  Kf(h)=\int_G k(g)\pi_g\,dg\,f(h)=\int_G k(g)f(hg)\,dg.
\end{equation}
Thus one may consider the composition of convolutions $K_2K_1$
with two kernels $k_{1}(g)$,  $k_{2}(g)\in \FSpace{L}{1}(G)$:
\begin{eqnarray*}
  K_2K_1f(h)&=&\int_G k_2(g_2) \int_G k_1(g_1)
  f(h g_2 g_1)\,dg_1\,dg_2\\
  &=&\int_G \left( \int_G k_2(g_2)\, k_1(g_2^{-1}g) \,dg_2\right)
  f(hg)\,dg
\end{eqnarray*}
where $g=g_2 g_1$. Thereafter it again is a convolution with the
kernel 
\begin{equation}\label{eq:comp}
  k(g)=\int_G k_2(g_2)\, k_1(g_2^{-1}g) \,dg_2.
\end{equation}
An important role in the consideration is played by the unique
solvability of equation $g_1g_2=g$ with respect to $g_1$.  For an
algebraic structure weaker than a group, if the inverse $g^{-1}$ is
not defined one can try to define a composition of convolutions by
the formula:
\begin{equation}\label{eq:blin}
  (k_1*k_2)(h)=\relstack{\int\!\int}{g_1 g_2=h} k_1(g_1)\,
  k_2(g_2)\,dg_1\,dg_2.
\end{equation}
However such a definition generates problems with an understanding of
the set double integral~\eqref{eq:blin} is taken over and the
corresponding measure. These difficulties disappear for left
c-semigroups: for given $b$ and $a$ there is at most one such
$x=[a^{-1}b]$ that $ax=b$. Thus we can preserve for c-semigroup $C$
definition~\eqref{eq:comp} with a small modification:
\begin{equation}\label{eq:cconv} 
  (k_2*k_1)(b)=\int_C   k_2(a)\, k_1([a^{-1}b])\,da 
\end{equation} 

To avoid problems with the definition of convolution we \emph{make
  an agreement} that for any function $f(c)$
\begin{equation}\label{eq:agreement}
  f([a^{-1}b])=0 \quad \textrm{ when }
  [a^{-1}b] \textrm{ does not exist in } C
\end{equation}
(i.e., the equation $ax=b$ does not have a solution). Our agreement is
equivalent to introduction of incidence algebra of functions $f(a,b)$,
such that $f(a,b)=0$ for $a>b$ for a poset (see Example~\ref{ex:poset}).
\begin{rem}
  As soon as one passes directly to convolution algebras and abandon
  points of c-semigroups themselves the future generalizations in the
  spirit of non-commutative geometry~\cite{Connes94} are possible.
\end{rem}
We give several examples of c-semigroup now.
\begin{example}[Principal]
  The main source of c-semigroups is the following construction. Let we
  have a set $S$ with an additional structure $\mathcal{A}$. Then the
  set of all mappings from $S$ \emph{in}to $S$, which preserve
  $\mathcal{A}$, forms a right c-semigroup.
\end{example}
\begin{example}
  If one considers in the previous Example only mappings $S$
  \emph{on}to $S$ then they give us a definition of groups.  Moreover
  a large set of c-semigroups can be constructed from groups directly.
  Namely, let $C$ be a subset of a group $G$ such that $(C\cdot
  C)\subset C$. Then $C$ with the multiplication induced from $G$ will
  be a c-semigroup. In such a way we obtain the c-semigroup of
  positive real (natural, rational) numbers with the usual
  addition. Entire (positive entire) numbers with multiplication also
  form a c-semigroup.  If $C$ contains the identity $e$ of $G$ then
  $e$ is a source of $C$.  Particularly, entire group $G$ is also a
  c-semigroup. We should note, that even for groups some of our
  technique will be new (for example \emph{tokens}), at least up to
  the author knowledge.
\end{example}

The last example demonstrates that we have many c-semigroups
with (left, right) invariant  measures such that Fubini's theorem
holds. This follows from corresponding constructions for groups.

We sign out the most important for the present paper example from the
described family and his alternative realization.
\begin{prop}\label{pr:series}
  The algebra of formal power series in a variable $t$ is
  topologically isomorphic to the convolution algebra
  $\FSpace{C}{}(\Space[+]{N}{})$ of continuous functions with
  point-wise convergence over the c-semigroup of non-negative integers
  $\Space[+]{N}{}$ with the discrete topology. The subalgebra of
  polynomials in one variable $t$ is isomorphic to convolution algebra
  $\FSpace{C}{0}(\Space[+]{N}{})$, which is dense in
  $\FSpace{C}{}(\Space[+]{N}{})$, of compactly supported functions on
  $\Space[+]{N}{}$.
\end{prop}
The proof is hardly needed. Note that
``puzzling questions'' about convergence of formal power series and
continuity of functionals on them correspond to the natural topology
on $\FSpace{C}{}(\Space[+]{N}{})$.

Less obvious family of examples can be constructed from a structure
quite remote (at least at the first glance) from groups.
\begin{example}\label{ex:poset}
  Let $P$ be a \emph{poset} (i.e., \emph{p}artialy \emph{o}rdered
  \emph{set}) and let $C$ denote the subset of Cartesian square
  $P\times P$, such that $(a,b)\in C$ iff $a\leq b$, $a,b\in P$. We
  can define a multiplication on $C$ by the formula:
\begin{equation}
  (a,b)(c,d)=\left\{\begin{array}{ll}
      \textrm{undefined }, & b \not = c;\\
      (a,d), & b=c.
    \end{array} \right.
\end{equation}
One can see that $C$ is a c-set. If $P$ is locally
finite, i.e., for any $a\leq b$, $a$, $b\in P$ the number of $z$ between
$a$ and $b$ ($a\leq z\leq b$) is finite, then we can define a
measure $d(a,b)=1$ on $C$ for any $(a,b)\in C$. With such a measure
\eqref{eq:cconv} defines the correct convolution on $C$:
\begin{equation}
  h(a,b)=\int_C f(c,d)\, g([(c,d)^{-1}(a,b)])\, d(c,d)
  =\sum_{a\leq z\leq b}  f(a,z)\, g(z,b).
\end{equation}
The constructed algebra is the fundamental \emph{incidence algebra} in
combinatorics~\cite{FoundationVI}.

Our present technique can be successfully applied only to
c-semigroups, not to c-sets. Thus it is an important observation
that \emph{reduced incidence algebra} construction~\cite[\S
~4]{FoundationVI} contracts c-sets to c-semigroups in many
important cases (however not in general). Particularly \emph{all
  combinatorial} reduced incidence algebras listed in
Examples~4.5--4.9 of~\cite{FoundationVI} are convolution algebras over
c-semigroups. However the incidence coefficients should be better
understood.
\end{example}

\begin{example}
  Any \emph{groupoid} (see~\cite{Weinstein96} and references herein)
  is a c-set. 
\end{example}

After given examples, it is reasonable to expect an applicability of
c-semigroups to combinatorics. We describe applications in the next
Section explicitly.

Now we point out some basic properties of c-semigroups, trying to be as
near as possible to their prototype---groups. All results will be
states for left c-semigroups. Under agreement~\eqref{eq:agreement} we
have
\begin{lem}[Shift invariance of integrals]
  \label{le:invariance}
  For any function $f(c)$ we have
\begin{equation}\label{eq:invariance}
  \int_C f(c)\, dc = \int_C f([a^{-1}b])\,db, \qquad \forall a\in C.
\end{equation}
\end{lem}
Proof follows from the observation that the integrals in both sides
of~\eqref{eq:invariance} is taken over the whole $C$ and measure $db$
is invariant. One should be warned that
\begin{displaymath}
  \int_C f(c)\, dc \neq \int_C f(ab)\,db
\end{displaymath}
for c-semigroups (unlike for groups) in general as well as
Lemma~\ref{le:invariance} is false for c-sets. Another very
useful property (which c-semigroups possess from the groups) is a
connection between linear functional and shift invariant operators.
Again unlike in the group case we should make the right choice of
shifts.
\begin{defn}
  We define \emph{left shift operator} $\lambda_{a}$ on the space of
  function $\FSpace{L}{1}(C)$ by the formula
  $[\lambda_{a}f](b)=f(ab)$.  A linear operator $L$ is (left)
  \emph{shift invariant} if $S \lambda_a =\lambda_a S $ for all $a\in
  C$. The \emph{augmentation} $\epsilon$ associated to a right
  source $e$ is the linear functional defined by $\epsilon(f)=f(e)$.
\end{defn}
\begin{prop}\label{pr:isomorphism}
  Let $C$ will be a c-semigroup with a right source $e$. The following
  three spaces are in one-to-one correspondence:
  \begin{enumerate}
  \item The space $\bos{L}$ of liner functionals on the space of
    functions $\FSpace{L}{}(C)$;
  \item The space $\bos{S}{}$ of shift invariant operators;
  \item The algebra of convolutions $\algebra{C}$ with functions from
    $\FSpace{L}{}^*(C)$. 
  \end{enumerate}    
  The correspondences are given by the formulas ($l\in
  \FSpace{L}{}^*(C)$, $S\in \FSpace{S}{}$, $k\in\algebra{C}$, $f\in
  \FSpace{L}{}(C)$)
  \begin{eqnarray*}
    l &\rightarrow& S, \mbox{ where } [Sf](a)=\scalar{l}{\lambda_a f}\\
    S &\rightarrow& l,\, \mbox{ where } lf=\epsilon(Sf)\\
    k &\rightarrow& l,\, \mbox{ where } lf=\int_C k(c) f(c)\, dc\\
    k &\rightarrow& S, \mbox{ where } Sf(a)=\int_C  k(c) f(ac)\, dc
  \end{eqnarray*}
\end{prop}
\begin{proof}
  Let we have a shift invariant operator $S$. Clearly that its
  composition $lf=\epsilon S f=[Sf](e)$ with the augmentation is a
  linear functional. By the definition of
  distributions~\cite[\S~III.4]{KirGvi82} such a linear functional can
  be represented as the integral $lk=\int_C k(c) f(c)\, dc$ with a
  (probably generalized) function $k(c)$ from the dual space. For any
  function $f$ we have $f(a)=\lambda_{ae^{-1}} f(e)$ thus by the
  shift-invariance of $S$
  \begin{displaymath}
    Sf(a)=[\lambda_{ae^{-1}} S f](e)= [S \lambda_{ae^{-1}} f](e)    
  \end{displaymath} 
  or by the definition of $l$ we find
  $Sf(a)=\scalar{l}{\lambda_{ae^{-1}} f}$. The integral form of the
  last identity is $Sk(a)=\int_C k(c) f(a c)\, dc$, which coincides
  with~\eqref{eq:conv}.
\end{proof} 
An essential r\^ole in all three approaches (operator,
Hopf algebras, semantic) to the umbral calculus is played by linear
functional (and thus associated linear shift invariant
operators). Connections between linear functionals and shifts can be
greatly simplified if we are able to express shifts via some ``linear
coefficients'' of functions. The known tool to do that is the Fourier
transform. An alternative tool for c-semigroups is presented in next
section.

\subsection{Tokens: Definition and Examples}
We formalize identity~\eqref{eq:binom}
\begin{displaymath}
  p_n(x+y)=\sum_{k=0}^{n} p_k(x) p_{n-k}(y).  
\end{displaymath} 
in the following notion, which will be useful in our
consideration and seems to be of an interest for a representation
theory. The main usage of this notion is $t$-transform defined in the
next subsection.
\begin{defn}
  Let $C_1$ and $C_2$ be two c-semigroups. We will say that a function
  $t(c_1,c_2)$ on $C_1 \times C_2$ is a \emph{token} from $C_1$ to
  $C_2$ if for any $c'_1 \in C_1$ and any $c_2, c_2' \in C_2$ we have
  \begin{equation}\label{eq:token}
    \int_{C_1} t(c_1,c_2)\, t([c_1^{-1}c_1'], c_2') \, dc_1 =
    t(c_1', c_2 c_2').
  \end{equation}
\end{defn}
\begin{rem}\label{re:fourier}
  We derive this definition ``experimentally'' from~\eqref{eq:binom} 
  (see also Examples bellow). However one can discover it
  ``theoretically'' as well. We already know from
  Proposition~\ref{pr:isomorphism} that a shift invariant operator $S$
  can be represented as the combination of shifts and a linear
  functional: $[Sf](a)=\scalar{l}{\lambda_a f}$.  This suggests to
  represent a function $f(b)$ as a linear combinations
  \begin{equation}\label{eq:linearcomb}
    f(b)=\int e(b,c) \, \hat{f}(c)\, dc
  \end{equation}
  of elementary components $e(b,c)$, which behave simply under shifts.
  Here ``simply'' means \emph{do not destroy linear
    combination}~\eqref{eq:linearcomb}.
  
  A way to achieve this is provided by the Fourier transform.  Here
  $e(b,c)=e^{ibc}$ and the action of a shift $\lambda_a$ reduces to
  multiplication:
  \begin{equation}\label{eq:exp}
    \lambda_a e(b,c)=e(b+a,c)=e^{i a c}e(b,c).
  \end{equation}
  This probably is the simplest solution but not the only possible one.
  More general transformation of this sort is given by the abstract
  wavelets (or coherent states):
  \begin{equation}\label{eq:wave}
    \lambda_a e(b,c)=e(b+a,c)=\int k(a,c,c_1)\, e(b,c_1) \,dc_1,
  \end{equation} 
  where we can particularly select
  $k(a,c,c_1)=e^{iac}\delta(c-c_1)$ to get the rule of
  exponents~\eqref{eq:exp} and $k(a,c,c_1)=e(c_1^{-1}c,a)$ to get the
  token property~\eqref{eq:token}.  Then shifts act on the functions
  defined via~\eqref{eq:linearcomb} as follows
  \begin{eqnarray*}
    \lambda_a f(b)&=&\lambda_a \int e(b,c)\, \hat{f}(c)\, dc\\
    &=&\int\!\!\int k(a,c,c_1)\,e(b,c_1) \,dc_1\, \hat{f}(c)\, dc\\
    &=&\int e(b,c_1) \left(\int k(a,c,c_1)\, \hat{f}(c)\,dc\right)  
    \,dc_1.
  \end{eqnarray*}
  Thus it again looks like~\eqref{eq:linearcomb} where the action of a
  shift reduced to an integral operator on the symbol $\hat{f}(c)$.
  In the case of token shifts affect in the way~\eqref{eq:token}, which
  stand between the extreme simplicity of~\eqref{eq:exp} and the almost
  unaccessible generality of~\eqref{eq:wave}.  This analogy with the
  Fourier transform will be employed in the next Subsection for the
  definition of $t$-transform. 
\end{rem}
\begin{rem}
  Relationships between harmonic analysis (group characters) and
  token-like structures (polynomial of binomial type) was already
  pointed in~\cite{BarBriRota80}.  Note that tokens are a complementary
  (in some sense) tool for the exponent $e^{ibc}$.  While both are
  particularly useful in an investigation of shift invariant operators
  they work in different ways. The Fourier transform maps shift
  invariant operators to operators of multiplication. $t$-Transform
  defined in the next subsection by means of tokens maps shift invariant
  operators on one c-semigroup to shift invariant operators on another
  c-semigroup.
\end{rem}
We present examples of tokens within classical objects.
\begin{example}\label{ex:polin1}
  Let $C_1$ be the
  c-semigroup of non-negative natural numbers $\Space[+]{N}{}$ and $C_2$
  be the c-semigroup (a group in fact) of real numbers both equipped by
  addition. Let $t(n,x)=p_n(x)$, where $p_n(x)$ is a polynomial
  sequence of binomial type. Then the characteristic property of such a
  sequence~\eqref{eq:binom} is equal to the definition of
  token~\eqref{eq:token} read from the right to left.
\end{example}
\begin{example}\label{ex:oper}
  We exchange the $C_1$ and $C_2$ from the previous Example, they are
  real numbers and non-negative natural numbers respectively. Let us
  take a shift invariant operator $L$. Any its power $L^n$ is again a
  shift invariant operator represented by a convolution with a
  function $l(n,x)$ with respect to variable $x$. The identity
  $L^{n+k}=L^{n} L^{k}$ can be expressed in the term of the
  correspondent kernels of convolutions:
  \begin{equation}
    \int_{\Space{R}{}} l(n,x)\, l(k,y-x)\, dx=l(n+k,y).
  \end{equation}
\end{example} 
First two examples with their explicit duality generate
a hope that our development can be of some interest for
combinatorics. But before working out this direction we would like to
present some examples of tokens outside combinatorics.
\begin{example}\label{ex:Poisson}
  Let $C_1$ be $\Space{R}{n}$ and $C_2$ is 
  $\Space{R}{n}\times\Space[+]{R}{}$---the ``upper half space'' in
  $\Space{R}{n+1}$. 
  For the space of harmonic function in
  $C_2=\Space{R}{n}\times\Space[+]{R}{}$  there is an integral
  representation over the boundary $C_1=\Space{R}{n}$:
\begin{displaymath}
  f(v,t)= \int_{C_1} P(u;v,t) f(u) \, du, \qquad u\in C_1, \ (v,t)\in
  C_2,\ v\in\Space{R}{n},\ t\in\Space[+]{R}{}.
\end{displaymath}
Here $P(u,v)$ is the celebrated Poisson kernel
\begin{displaymath}
  P(u;v,t)=\frac{2}{\modulus{S_n}}\frac{t}{(\modulus{u-v}^2+t^2)^{(n+1
      )/2}}
\end{displaymath}
with the property usually referred as a \emph{semigroup
  property}~\cite[Chap.~3, Prob.~1]{Akhiezer88}
\begin{displaymath}
  P(u;v+v', t+t')=\int_{C_1} P(u';v',t') P(u-u';v,t)\, du'.
\end{displaymath}
We meet the token in analysis.
\end{example}
\begin{example}\label{ex:Weierstrass}
  We preserve the meaning of $C_1$ and $C_2$ from the previous Example
  and define the Weierstrass (or Gauss-Weierstrass) kernel by the
  formula:
  \begin{displaymath}
     W(z;w,\tau)=\frac{1}{(\sqrt{2\pi \tau})^n} e^{\frac{-
        \modulus{z-w}^2}{2\tau}}=
    \frac{1}{(2\pi)^n}\int_{\Space{R}{n} }
    e^{-\frac{\tau}{2}\modulus{u}^2} e^{-(u,z-w)}\, du,
  \end{displaymath}
  where $z\in C_1$, $(w,\tau)\in C_2$. Function $W(z;w,\tau)$ is the
  fundamental solution to the heat equation~\cite[\S~2.3]{Evans98}.
  We again have~\cite[Chap.~3, Prob.~1]{Akhiezer88}
  \begin{displaymath}
    W(z;w+w',\tau +\tau')=\int_{C_1} W(z';w',\tau') W(z-z';w,\tau)\, dz.
  \end{displaymath}
  Thus we again meet a token.
\end{example} 
Two last Examples open a huge list of integral
kernels~\cite{Akhiezer88} which are tokens of
analysis. We will see the r\^ole of tokens for analytic function
theory and reproducing kernels later.

We finish the subsection by the following elementary property of
token.
\begin{lem}\label{le:t-reprod}
  Let c-semigroup $C_2$ has a left source $e_2$. Then
  $k(c_1)=t(c_1,e_2)$ has the reproducing property
  \begin{displaymath}
    t(c_1',c_2)=\int_{C_1} k(c_1)\,
    t([c_1^{-1}c_1'],[e_2^{-1}c_2])\,dc_1 
    \quad    \forall c'_1\in C_1,\  c_2\in C_2.
  \end{displaymath}
  Particularly if $e_2\in C_2$ is a left unit then:
  \begin{displaymath}
    t(c_1',c_2)=\int_{C_1} k(c_1)\,
    t([c_1^{-1}c_1'],c_2)\,dc_1 
    \quad    \forall c'_1\in C_1,\ c_2\in C_2.
  \end{displaymath}
\end{lem}
\begin{proof}
  We apply the characteristic property of a token~\eqref{eq:token} for
  particularly selected elements in $C$:
  \begin{eqnarray*}
    t(c_1',c_2)&=&t(c_1',e_2[e_2^{-1}c_2])\\
    &=&\int_{C_1} t(c_1,e_2) t([c_1^{-1}c_1'],[e_2^{-1}c_2])\,dc_1\\
    &=&\int_{C_1} k(c_1) t([c_1^{-1}c_1'],[e_2^{-1}c_2])\,dc_1.
  \end{eqnarray*}
\end{proof}

\subsection{t-Transform}

We would like now to introduce a transformation associated with tokens
(t-transform). t-Transform is similar to the Fourier and wavelet
transforms. The main difference---we do not insist that our
transformation maps convolutions to multiplications. We will be much
more modest---our transformation maps convolutions to other
convolutions, eventually more simple or appropriate.

\comment{\begin{eqnarray} [K t](c_1',c_2) &=&\int_{C_2} k(c_2')
    t(c_1',c_2 c_2') \, dc_2'
    \nonumber \\
    &=& \int_{C_2} k(c_2') \int_{C_1} t(c_1,c_2) t([c_1^{-1}c_1'],
    c_2') \, dc_1
    \, dc_2'\quad \  \label{eq:transform01}  \\
    &=& \int_{C_1} \int_{C_2} k(c_2') t([c_1^{-1}c_1'], c_2') \, dc_2'
    t(c_1,c_2) \, dc_1 \quad \  \nonumber \\
    &=& \int_{C_1} \hat{k}([c_1^{-1}c_1']) t(c_1,c_2) \,dc_1
    \label{eq:transform03}
\end{eqnarray}
where
\begin{displaymath}
  \hat{k}(c_1) = \int_{C_2} k(c_2')\, t(c_1, c_2') \, dc_2'
\end{displaymath}}\label{comment:one}

Having a function $\hat{f}(c_1)$ on $C_1$ we can consider it as
coefficients (or a \emph{symbol}) and define the function $f(c_2)$ on
$C_2$ by means of the token $t(c_1,c_2)$:
\begin{equation}\label{eq:tf-transform}
  f(c_2)=\int_{C_1}  \hat{f}(c_1)\, t(c_1,c_2)\, dc_1.
\end{equation} 

We will say that $f(c_2)\in\FSpace{A}{}(C_2)$ if in~\eqref{eq:tf-transform}
$\hat{f}(c_1)\in\FSpace{C}{0}(C_1)$ the space of continuous functions
with a compact support. So
we can consider the token $t(c_1,c_2)$ as a kernel of integral
transform $\FSpace{C}{0}(C_1) \rightarrow \FSpace{A}{}(C_2)$. The use
of the notation reserved for the Fourier transform
(\,$\hat{}$\,---``hat'') is justified (at least partially) by
Remark~\ref{re:fourier}. We would like to mention the following
\begin{cor}\label{co:reprod}
  The reproducing property of $k(c_1)$ from Lemma~\ref{le:t-reprod}
  can be extended by linearity to $\FSpace{A}{}(C_1)$:
  \begin{displaymath}
    f(c_1')=\int_{C_1} k(c_1)\, f([c_1^{-1}c_1'])\,dc_1.
  \end{displaymath}
\end{cor}

\begin{example}\label{ex:polin2}
  We continuing with a polynomial sequence of binomial type $p_k(x)$,
  which was presented as a token $p(k,x)$ in Example~\ref{ex:polin1}.
  Having few numbers $\hat{a}_1$, $\hat{a}_2$, \ldots,
  which are interpreted as a function
  $\hat{a}(k)$ in $\FSpace{C}{0}(\Space[+]{N}{})$, we can
  construct a polynomial
  \begin{displaymath}
    a(x)=\int_{\Space[+]{N}{}} \hat{a}(k)\,p(k,x)\,dk=\sum_k
    \hat{a}_k p_k(x),
  \end{displaymath}
  that is a function on $\Space{R}{}$. 
\end{example}
\begin{example}
  Let $\delta^{(k)}(x)$ be a token from $\Space{R}{}$ to
  $\Space[+]{N}{}$ (Example~\ref{ex:oper}). Then for a smooth at point 
  $0$ function $f(x)$ the integral transformation 
  \begin{displaymath}
    f^{(k)}(0)=\int_{\Space{R}{}} f(x)\, \delta^{(k)}(x)\, dx
  \end{displaymath}
  produces a function on $\Space[+]{N}{}$.
\end{example}
\begin{example}
  Let a real valued function $f(u)$ is defined on
  $C_1=\Space{R}{n}$. Then we can define the function $Pf(v,t)$ on
  $C_2=\Space{R}{n}\times\Space[+]{R}{}$ by means of the Poisson
  integral:
  \begin{eqnarray*}
    Pf(v,t) & = & \int_{C_1} P(u; v,t) f(u)\, du\\
    &=& \frac{2}{\modulus{S_n}}\int_{\Space{R}{n}}  
    \frac{t}{(\modulus{u-v}^2+t^2)^{(n+1)/2}} f(u)\,du. 
  \end{eqnarray*}
  The image is a harmonic function solving the Dirichlet
  problem~\cite[Appendix, \S~2]{ShabatI}, \cite[\S~2.2.4]{Evans98}.
\end{example}
\begin{example}
  Similarly we define the Weierstrass transform:
  \begin{eqnarray*}
    Wf(w,\tau)&=&\int_{C_1}  W(z;w,\tau) f(z)\,dz\\
    &=& \frac{1}{(\sqrt{2\pi \tau})^n}\int_{\Space{R}{n}}  e^{\frac{-
        \modulus{z-w}^2}{2\tau}} f(z)\,dz
  \end{eqnarray*}
  The image is a function satisfying the heat equation.
\end{example}
Let $K$ be a convolution on $C_2$ with a kernel $k(c_2)$. Then it acts
on function $f(c_2)$ defined by~\eqref{eq:tf-transform} as follows
\begin{eqnarray}
  [Kf](c_2) &=& \int_{C_2} k(c_2') f(c_2 c_2')\, dc_2' \nonumber \\
  &=& \int_{C_2} k(c_2') \int_{C_1} \hat{f}(c_1') t(c_1',c_2 c_2')
  \, dc_1' \, dc_2' \nonumber \\
  &=& \int_{C_2} k(c_2') \int_{C_1}  \hat{f}(c_1')
  \int_{C_1} t(c_1,c_2) t([c_1^{-1}c_1'], c_2') \, dc_1
  \, dc_1'\, dc_2'\qquad \  \label{eq:transform11}  \\
  &=& \int_{C_1} \int_{C_1} \int_{C_2} k(c_2')\, t([c_1^{-1}c_1'], c_2')
  \, dc_2'\, \hat{f}(c_1') \, dc_1'\, t(c_1,c_2)
  \,dc_1  \nonumber \\
  &=& \int_{C_1} \int_{C_1}  \int_{C_2} k(c_2')\, t(c_1'', c_2') \, dc_2'\,
  \hat{f}(c_1 c_1'') \, dc_1''\, t(c_1,c_2)
  \,dc_1  \label{eq:transform12} \\
  &=& \int_{C_1} \left(\int_{C_1}  [\oper{T}k](c_1'')
    \hat{f}(c_1 c_1'') \, dc_1'' \,\right) t(c_1,c_2)
  \,dc_1  \label{eq:transform13}
\end{eqnarray}
where
\begin{equation}\label{eq:kernel}
  [\oper{T}k](c_1'')=\int_{C_2} k(c_2')\, t(c_1'', c_2') \, dc_2'.
\end{equation}
We use in~\eqref{eq:transform11} the characteristic property of the
tokens~\eqref{eq:token} and the change of variables
$c_1''=[c_1^{-1}c_1']$ in~\eqref{eq:transform12}.  Note that big
brasses in~\eqref{eq:transform13} contain a convolution on $C_1$ with
the kernel~\eqref{eq:kernel}.  So
transformation~\eqref{eq:tf-transform} maps convolutions on $C_2$ to
the convolutions on $C_1$. The mapping deserves a special name.

\begin{defn}
  t-Transform is a linear integral transformation $\algebra{C}(C_2)
  \rightarrow \algebra{C}(C_1)$ defined on kernels by the formula:
\begin{equation}\label{eq:t-transform}
  [\oper{T}k](c_1)=\int_{C_2} t(c_1,c_2)\, k(c_2)\, dc_2.
\end{equation}
\end{defn}
According to Lemma~\ref{le:invariance} we have
\begin{equation}\label{eq:token-inv}
  [\oper{T}k](c_1)=\int_{C_2} t(c_1,[c^{-1}_2 c'_2])\, k([c_2^{-1} c'_2])\,
  dc'_2, \qquad \forall c_2\in C_2.
\end{equation}
The t-transform is more than just a linear map:
\begin{thm}\label{th:t-transform}
  t-transform is an algebra homomorphism of the convolution algebra
  $\algebra{C}(C_2)$ to the convolution algebra $\algebra{C}(C_1)$,
  namely
  \begin{equation}\label{eq:t-conv}
    [\oper{T}(k_1*k_2)]=[\oper{T}k_1]*[\oper{T}k_2],
  \end{equation} 
  where $*$ in the left-hand side denotes the
  composition of convolutions on $C_2$ and in the right-hand side---on
  $C_1$.
\end{thm}
\begin{proof}
  We silently assume that conditions of Fubini's theorem fulfill. It
  turns to be that Fubini's theorem will be our main tool, which is used
  in almost all proofs. Thereafter we assume that it fulfills through
  entire paper. We have:
  \begin{eqnarray}
    [\oper{T}(k_2*k_1)](c'_1)&=&\int_{C_2} t(c'_1,c'_2) \int_{C_2}
    k_2(c_2)\, k_1([c_2^{-1}c'_2]) \,dc_2\,dc'_2\nonumber \\
    &=&\int_{C_2} \int_{C_2} \int_{C_1}  t(c_1,c_2)\, 
      t([c_1^{-1}c'_1],[c_2^{-1}c'_2])\, dc_1 \label{eq:substit1}\\
    && \qquad\times k_2(c_2)\, k_1([c_2^{-1}c'_2])\,dc_2\,dc'_2 \nonumber\\
    &=& \int_{C_1} \int_{C_2}\, t(c_1,c_2)\, k_2(c_2) \,dc_2 \nonumber\\
    && \qquad \times \int_{C_2} t([c_1^{-1}c'_1],[c_2^{-1}c'_2])\,  
    k_1([c_2^{-1}c'_2])\,dc'_2\, dc_1 \nonumber \\
    &=& \int_{C_1} [\oper{T}k_2](c_1) \, [\oper{T}k_1]([c_1^{-1}c'_1])
    \,dc_1\label{eq:substit2}.
  \end{eqnarray}
  Of course, we use in transformation~\eqref{eq:substit1} the
  characteristic property of the tokens~\eqref{eq:token} and
  conclusion~\eqref{eq:substit2} is based on~\eqref{eq:token-inv}.
\end{proof}
Knowing the connection between shift invariant operators and
convolutions (Proposition~\ref{pr:isomorphism}) we derive
\begin{cor}\label{co:t-shift}
  t-Transform maps left shift invariant operators on $C_2$ to left
  shift invariant operators on $C_1$.
\end{cor}

\section{Shift Invariant Operators}\label{sec:shift-invar-oper}
In this Section we present
some results in terms of c-semigroups, tokens and $t$-transform
and illustrate by examples from combinatorics and analysis.

\subsection{Delta Families and Basic Distributions}
This Subsection follows the paper~\cite{KahOdlRota73} with appropriate
modifications.
\begin{defn}
  Let $C_1$ and $C_2$ be two c-semigroups; let $e_1$ be a right source
  of $C_1$. We will say that linear shift invariant operators $Q(c_1)$,
  $c_1\in C_1$ on $\FSpace{A}{}(C_2)$ over c-semigroup $C_2$ form a
  \emph{delta family} if 
  \begin{enumerate}
  \item for a token $t(c_1,c_2)$
    \begin{displaymath}
      Q(c_1)\, t(c_1, c_2)= t(e_1, c_2) \quad 
      \textrm{ for all } \quad c_2\in C_2.
    \end{displaymath}
  \item the semigroup property $Q(c_1)\, Q(c_1')=Q(c_1c_1')$ holds.
  \end{enumerate}
\end{defn}
\begin{example}\label{eq:polin4}
  Let $C_1$ and $C_2$ be as in Examples~\ref{ex:polin1},
  \ref{ex:polin2}, and \ref{ex:polin3}. Let also
  $t(n,x)=\frac{x^n}{n!}$. Then the family of positive integer powers
  of the derivative $Q(n)=D^n$ forms a delta. Other deltas for the
  given $C_1$ and $C_2$ are listed in~\cite[\S~3]{FoundationIII}. They
  are positive integer powers of difference, backward difference,
  central difference, Laguerre, and Abel operators.
\end{example}
\begin{example}\label{ex:Poisson-delta}
  For the Poisson kernel from Example~\ref{ex:Poisson} we can
  introduce delta family as Euclidean shift operators: 
  \begin{eqnordisp}
    Q(x)\, f(v,t)= f(v+x,t), \qquad x\in \Space{R}{n}.
  \end{eqnordisp}
\end{example}
\begin{prop}\label{pr:kernel-df}
  Let $q(c_1,c_2)$ be the kernel corresponding to a delta family
  $Q(c_1)$ as a convolution over $C_2$. Then $q(c_1,c_2)$ is a token
  from $C_2$ to $C_1$. Namely
\begin{displaymath}
  q(c_1c_1',c_2')=\int_{C_2} q(c_1,c_2)\, q(c_1',[c_2^{-1}c_2'])\,dc_2.
\end{displaymath}
\end{prop}
\begin{proof}
  This is a simple restatement in terms of kernels of the semigroup
  property $Q(c_1)\,Q(c_1')=Q(c_1 c_1')$.
\end{proof}
\begin{defn}
  The function $t(c_1,c_2)$ is called the \emph{basic distribution}
  for a delta family $Q(c_1)$ if:
\begin{enumerate}
\item $t(c_1,e_2)= k(c_1)$ has the reproducing property from
  Lemma~\ref{le:t-reprod} over $C_1$ (particularly if 
  $k(c_1)=\delta(c_1)$);
\item $Q(c_1)\, t(c_1',c_2)=t([c_1^{-1}c_1'],c_2)$.
\end{enumerate}
\end{defn}
\begin{example}
  For the classic umbral calculus this conditions turn to be
  $p_0(x)=1$; $p_n(0)=0$, $n>0$; $Qp_n(x)=p_{n-1}(x)$.
\end{example}
\begin{example}
  For the Poisson integral the reproducing property can be obtained
  if we consider a generalized Hardy space
  $\FSpace{H}{2}(\Space[+]{R}{n+1})$ of $\FSpace{L}{2}$-integrable
  functions on $\Space{R}{n}$, which are limit value of harmonic
  functions in the upper half space
  $\Space{R}{n}\times\Space[+]{R}{}$. Then kernel $k(u-u')$ is defined
  as a distribution given by the limit of integral on
  $\FSpace{H}{2}(\Space{R}{n})$:
  \begin{eqnarray*}
    f(u) & = & \int_{\Space{R}{n}} k(u-u')\, f(u')\,du'\\
    & =& \frac{2}{\modulus{S_n}} \,
    \lim_{(v,t) \rightarrow (0,0)} \,\int_{\Space{R}{n}}  
    \frac{t}{(\modulus{(u-u')-v}^2+t^2)^{(n+1)/2}} f(u')\,du'. 
  \end{eqnarray*}
  This follows from the boundary property of the Poisson
  integral~\cite[Appendix, \S~2]{ShabatI}, \cite[\S~2.2.4]{Evans98}. 
\end{example}
\begin{thm}\label{th:b-dist}
  Let $C_1$ and $C_2$ be c-semigroups.
  \begin{enumerate}
  \item If $t(c_1,c_2)$ is a basic distribution for some delta family
    $Q(c_1)$, $c_1\in C_1$, then $t(c_1,c_2)$ is a token from $C_1$ to
    $C_2$.
  \item If $t(c_1,c_2)$ is a token from $C_1$ to $C_2$, then it is a
    basic distribution for some delta family $Q(c_1)$, $c_1\in C_1$.
  \end{enumerate}
\end{thm}
\begin{proof}
  Let $t(c_1,c_2)$ be a basic distribution for $Q(c_1)$ and let
  $q(c_2,c_1)$ be a kernel of operator $Q(c_1)$ as a convolution on
  $C_2$. This means that
\begin{displaymath}
  \scalar[c_2]{q(c_2,c_1)}{t(c_1',c_2)}=k ([c_1^{-1}c_1'])
\end{displaymath}
has a reproducing property due to Lemma~\ref{le:t-reprod} and
Corollary~\ref{co:reprod}.  Thus we can trivially express
$t(c'_1,c'_2)$ as the integral
\begin{displaymath}
  t(c'_1,c'_2)=\int_{C_1} t(c_1,c_2')
  \scalar[c_2]{q(c_2,c_1)}{t(c_1',c_2)}\, dc_1.
\end{displaymath}
By linearity we have a similar expression for any function in
$\FSpace{A}{}(C_2)$:
\begin{displaymath}
  f(c'_2)=\int_{C_1} t(c_1,c_2') \scalar[c_2]{q(c_2,c_1)}{f(c_2)}\,
  dc_1.
\end{displaymath}
In particular
\begin{displaymath}
  t(c_1',c_2''c_2')=\int_{C_1} t(c_1,c_2')
  \scalar[c_2]{q(c_2,c_1)}{t(c_1',c_2''c_2)}\,
  dc_1.
\end{displaymath}
But
\begin{eqnarray*}
  \scalar[c_2]{q(c_2,c_1)}{t(c_1',c_2''c_2)}
  &=&\intersect{Q(c_1)\lambda_{c_2''}t(c_1',c_2)}{c_2=e_2}\\
  &=&\intersect{\lambda_{c_2''}Q(c_1)t(c_1',c_2)}{c_2=e_2}\\
  &=&\intersect{\lambda_{c_2''}t([c_1^{-1}c_1'],c_2)}{c_2=e_2}\\
  &=&t([c_1^{-1}c_1'],c''_2),
\end{eqnarray*}
and therefore
\begin{displaymath}
  t(c_1',c_2''c_2')=\int_{C_1} t(c_1,c_2')\,
  t([c_1^{-1}c_1'],c''_2)\, dc_1;
\end{displaymath}
that is, the distribution $t(c_1,c_2)$ is a token.

Suppose now that $t(c_1,c_2)$ is a token. We define a family of
operators $Q(c_1)$ on $\FSpace{A}{}(C_2)$ by identities
\begin{displaymath}
  Q(c_1)\, t(c_1',c_2)=t([c_1^{-1}c_1'],c_2)
\end{displaymath}
and extending to whole $\FSpace{A}{}(C_2)$ by the linearity. The
semigroup property $Q(c_1)Q(c_2)=Q(c_1c_2)$ follows automatically. The
main point is to show that $Q(c_1)$ are shift invariant. We may
trivially rewrite the characteristic property of token as
\begin{displaymath}
  t(c_1',c_2'c_2)=\int_{C_1} t(c_1,c_2')\, Q(c_1)\, t(c_1',c_2)\, dc_1,
\end{displaymath}
which can be extended by linearity to any function in
$\FSpace{A}{}(C_2)$:
\begin{displaymath}
  f(c_2'c_2)=\int_{C_1} t(c_1,c_2')\, Q(c_1)\,f(c_2)\, dc_1,
\end{displaymath}
Now replace $f$ by $Q(c_1')f$
\begin{displaymath}
  Q(c_1')f(c_2'c_2)=\int_{C_1} t(c_1,c_2')\, Q(c_1'c_1)\,f(c_2)\, dc_1,
\end{displaymath}
But the left-hand side of the previous identity is nothing else as
$[\lambda_{c_2'}Q(c_1')f](c_2)$ and the right-hand side is
\begin{eqnarray*}
  \int_{C_1} t(c_1,c_2')\, Q(c_1'c_1)\,f(c_2)\, dc_1
  &=&Q(c_1')\left( \int_{C_1} t(c_1,c_2')\, Q(c_1)f(c_2)\, dc_1\right)\\
  &=&Q(c_1')\left( f(c_2' c_2)\right)= [Q(c_1')\lambda_{c_2'}f](c_2),
\end{eqnarray*}
i.e., $Q$ is shift invariant.
\end{proof}
The following Lemma is obvious.
\begin{lem}
  Let $q(c_2,c_1)$ be the kernel of $Q(c_1)$ as a convolution on
  $C_2$, then the t-transform for the token $q(c_2,c_1)$ is a
  right inverse operator for t-transform~\eqref{eq:t-transform} with
  respect to $t(c_1,c_2)$.
\end{lem}
The above Lemmas justify the following Definition.
\begin{defn}
  Let $t(c_1,c_2)$ be the basic distribution for a delta family
  $Q(c_1)$. Then we call $t(c_1,c_2)$, which is a token by
  Theorem~\ref{th:b-dist}, and the kernel $q(c_2,c_1)$ of $Q(c_1)$,
  which is a token by Proposition~\ref{pr:kernel-df}, \emph{dual
    tokens}.
\end{defn}
\begin{example}\label{ex:xn-diracs}
  The polynomials $t(n,x)=x^n/n!$ and distributions
  $q(n,x)=\delta^{(n)}(x)$ (which are kernels of derivatives
  operators $D^{(k)}$) are canonical dual tokens.  
\end{example}

As the reader may see in Theorem~\ref{th:b-dist} we almost do not
change proofs of~\cite[Theorem~1, \S~2]{KahOdlRota73}. So we give our
version of next two results without proofs.
\begin{thm}[First Expansion Theorem]
  Let $S$ be a shift invariant operator on $\FSpace{A}{}(C_2)$ and let
  $Q(c_1)$, $c_1\in C_1$ be a delta family with basic distribution
  $t(c_1,c_2)$. Then
  \begin{displaymath}
    S=\int_{C_1} a(c_1)\, Q(c_1)\,dc_1
  \end{displaymath}
  with
  \begin{displaymath}
    a(c_1)=\scalar[c_2]{s(c_2)}{t(c_1,c_2)}=\intersect{St(c_1,c_2)}{c_2=e_2}.
  \end{displaymath}
\end{thm}
\begin{thm}[Isomorphism Theorem]\label{th:isomorphism}
  Let $Q(c_1)$, $c_1\in C_1$ be a delta family with basic distribution
  $t(c_1,c_2)$. Then mapping which carries a shift invariant operator
  $S$ to a function
\begin{displaymath}
  s(c_2)=\int_{C_1} a(c_1)\, q(c_2,c_1)\,dc_1
\end{displaymath}
with
\begin{displaymath}
  a(c_1)=\scalar[c_2]{s(c_2)}{t(c_1,c_2)}=\intersect{St(c_1,c_2)}{c_2=e_2}.
\end{displaymath}
is an isomorphism of operator algebra to a convolution algebra on
$C_2$.
\end{thm}
\begin{example}
  For a sequence of binomial type $p_n(x)$ and the associated delta
  operator $Q$ we can decompose a shift invariant operator $S$ on the
  space of polynomials as
  \begin{eqnordisp}
    S= \sum_{k=0}^\infty a_k Q^k, \qquad \textrm{ where } \quad
    a_k=[Sp_k(x)]_{x=0}. 
  \end{eqnordisp}
\end{example}
\begin{example}
  Let a shift invariant operator $S$ on
  $\FSpace{H}{2}(\Space[+]{R}{n+1})$ is given by a kernel $s(v,t)$:
  \begin{eqnordisp}{}
    [S f](v,t)= \int_{\Space[+]{R}{n+1}} s(v',t')
    f(v+v',t+t')\,dv'\,dt'. 
  \end{eqnordisp} Then we can represent $S$ by means of the delta
  family $Q(u)$, $u\in\Space{R}{n}$ from
  Example~\ref{ex:Poisson-delta} as follows:
  \begin{eqnordisp}
    S = \int_{\Space{R}{n}} a(u)\, Q(u)\,du
  \end{eqnordisp}
  where
  \begin{eqnordisp}
    a(u)= \int_{\Space[+]{R}{n+1}} s(v,t)\, P(u;v,t)\, dv\,dt.
  \end{eqnordisp}
  And a similar representation is true for the Weierstrass kernel.
\end{example}

The following result has a somewhat cumbersome formulation and a
completely evident proof. In the classic umbral calculus case the
formulation turns to be very natural, but proof more hidden (see next
example).
\begin{thm}\label{th:order}
  Let $t_1(c_1,c_2)$ and $t_2(c_1,c_2)$ be two tokens from $C_1$ to
  $C_2$. Let c-semigroup $C_1$ be equipped with an partial order
  relation $>$.  Moreover let
  \begin{equation}\label{eq:order}
    t_2(c_1,c_2)=\int_{C_1} a(c_1,c_1')\, t_1(c_1',c_2)\,dc_1'
  \end{equation}
  where $a(c_1,c_1)=0$ whenever $c_1>c_1'$. Let $Q_2(c_1)$ be a delta
  family associated to $t_2(c_1,c_2)$. Let
  $\object{sup}_{t_1}(f)\in C_1$ denotes the support of
  $t$-transform of a function $f(c_2)$ with respect to $t_1(c_1,c_2)$.
  Then for any function $f(c_2)$ and $c_1\in C_1$:
  \begin{displaymath}
    \lambda_{c_1}\left(\object{sup}_{t_1}(Q(c_1)f)\right) 
    \subset \object{sup}_{t_1}(f).
  \end{displaymath}
\end{thm}
\begin{proof}
  This is a direct consequence of shift invariance of a delta family
  $Q(c_1)$ as operators on $C_2$ and Corollary~\ref{co:t-shift}.
\end{proof}
\begin{example}
  C-semigroup $\Space[+]{N}{}$ of natural numbers is naturally ordered.
  By the very definition~\cite[\S~3]{KahOdlRota73} for any polynomial
  sequence of binomial type $p(k,x)$ is exactly of degree $k$ for all
  $k$. Thus for any two such sequences $p_1(k,x)$ and $p_2(k',x)$ we have
  \begin{displaymath}
    p_1(k,x)=\sum_{k'=1}^k a(k,k')\, p_2(k',x),
  \end{displaymath} 
  i.e., \eqref{eq:order} is satisfied. Let now $Q(k)=Q^k$ be the delta
  family associated to $p_1(k,x)$, then by Theorem~\ref{th:order}
  polynomial $Q(k)p_2(k',x)$ if of degree $k'-k$. This
  gives~\cite[Propositions~1 and~2, \S~2]{KahOdlRota73}.
\end{example}
Generally speaking an order on a set of combinatorial number allows to
consider some recurrence relations, which express a combinatorial
number via smaller numbers of the same kind. We will return to this
subject in Subsection~\ref{ss:recurrence} in a connection with
generating functions.

\subsection{Generating Functions and Umbral 
  Functionals}\label{ss:generating} 

It is known that generating functions is a powerful tool in
combinatorics. 
\comment{It was pointed
in~\cite{FoundationVI}:
\begin{quote}
  \ldots too often generating functions have been considered to be
  simply an application of the current methods of harmonic analysis.
  From several of the examples discussed in this paper it will appear
  that this not the case \ldots
\end{quote} 
In the present paper we extend borders of \emph{the current methods of
  harmonic analysis} introducing the notion of tokens, which is
similar to a group character. } We try to interpret the notion of
generating functions in our terms.
Particularly we will treat Examples 4.5--4.9 from~\cite{FoundationVI}
accordingly.
\begin{defn}
  Let function $f(c_1)$ represents a set of combinatorial quantities
  indexed by points of c-semigroup $C_1$ (a \emph{combinatorial
    function} for short).  Let $C_2$ be another c-semigroup and
  $t(c_1,c_2)$ be a token between them. The t-transform
\begin{displaymath}
  \hat{f}(c_2)=\int_{C_1} f(c_1)\, t(c_1,c_2) \,dc_1
\end{displaymath}
of $f(c_1)$ is called \emph{generating function} for $f(c_1)$ (with
respect to $C_2$ and $t(c_1,c_2)$)
\end{defn}

Following~\cite{Rota64a} an exceedingly useful method for defining a
function $f(c_1)$ is to apply an umbral linear functional $l$ to a
token $p(c_1,c_2)$:
\begin{equation}\label{eq:umbra}
  f(c_1)=\scalar{l}{p(c_1,c_2)}:=\int_{C_2} l(c_2)\, p(c_1,c_2)\,dc_2.
\end{equation} 
If we now consider the generating function $\hat{f}(c_2)$ with
respect to the token $q(c_2,c_1)$, which is dual to $p(c_1,c_2)$ then
we will found that
\begin{eqnarray*}
  \hat{f}(c_2)&=&\int_{C_1} f(c_1)\, q(c_2,c_1)\,dc_1\\
  &=&\int_{C_1} \int_{C_2} l(c'_2)\,p(c_1,c'_2)\,dc'_2 \,
  q(c_2,c_1)\,dc_1\\
  &=&\int_{C_2} l(c'_2) \int_{C_1} p(c_1,c'_2)\, 
  q(c_2,c_1)\,dc_1\,dc'_2\\
  &=&\int_{C_2} l(c'_2)\, \delta (c'_2,c_2)\,dc'_2\\
  &=& l(c_2)
\end{eqnarray*}
So we found the following simple connection between umbral functionals
and generating functions
\begin{thm}
  The generating function $\hat{f}(c_2)$ for a combinatorial
  function $f(c_1)$ with respect to a token $q(c_2,c_1)$ is the kernel
  of umbral functional for $f(c_1)$ with respect to the token
  $p(c_1,c_2)$ dual to token $q(c_2,c_1)$.
\end{thm}
\begin{example}
  Let $\{f_n\}_{k=0}^\infty$ be the sequence of the Fibonacci
  numbers. Let $F(x)$ be the generating function associated to
  $\{f_n\}_{k=0}^\infty$ by the token $q(k,x)=\delta^{(k)}(x)$:
  \begin{eqnarray*}
    F(x) & = & \sum_{k=0}^{\infty} f_n\, \delta^{(k)}(x) \qquad 
    \left[ = \frac{1}{2\pi}
      \int_{-\infty}^{\infty} \frac{1}{1-(i\xi)-(i\xi)^2}\,
      e^{-i\xi x}\,d\xi\right]\\
     & = & \frac {1}{\sqrt{5}} \,e^{\frac{1-\sqrt{5}}{2}x} \,
  \left(1+ (e^{ \sqrt{5}\,x}-1)\,\chi(-x)\right)
  \end{eqnarray*}
  where $\chi(x)$ is the Heaviside function:
  \begin{eqnordisp}
    \chi(x)=\left\{
    \begin{array}{ll}
      0, &x<0; \\ 1, & x\geq 0.
    \end{array}
    \right.
  \end{eqnordisp}
  Then $F(x)$ is the kernel of an umbral functional, which generates the
  sequence $\{f_n\}_{k=0}^\infty$ from polynomial sequence
  $t(n,x)=x^n/n!$ which is dual to $q(k,x)=\delta^{(k)}(x)$ 
  (Example~\ref{ex:xn-diracs}):
  \begin{eqnordisp}
    f_n = \int_{-\infty}^\infty F(x)\, \frac{x^n}{n!}\, dx.
  \end{eqnordisp}
  In this way we obtain a constructive solution for the classic moment
  problem~\cite{Akhiezer61}. 
\end{example}

\subsection{Recurrence Operators and Generating 
  Functions}\label{ss:recurrence} 

We are looking for reasons why some generating functions are better
suited to handle given combinatorial functions than other.  The
answer is given by their connections with recurrence operators. 

We will speak on \emph{combinatorial functions} here that are
functions defined on a c-semigroups. Different numbers (Bell,
Fibonacci, Stirling, etc.) known in combinatorics are, of course,
combinatorial functions defined on c-semigroup $\Space{N}[+]{}$.
\begin{defn}
  An operator $R$ on $C_1$ is said to be a \emph{recurrence operator}
  with respect to a token $t(c_1,c_2)$ for a combinatorial function
  $f(c_1)$ if
  \begin{equation}\label{eq:recurrent}
    [Rf](c_1)=t(c_1,e_2).
  \end{equation}
  Particularly $[Rf](c_1)$ has a reproducing property from
  Lemma~\ref{le:t-reprod}.
\end{defn}

If we consider an algebra of convolutions of combinatorial functions
then by the general property of t-transform \emph{it is isomorphic to
  convolution algebra of their generating function}, which
occasionally can be isomorphic to an algebra with multiplication
(multiplication of formal power series, exponential power series,
Dirichlet series, etc.). This was already observed by other means
in~\cite{FoundationVI}.

The following Theorem connects generating functions and recurrent
operators.
\begin{thm}
  \label{th:wisdom-teeth}
  Let $C_1$ and $C_2$ be c-semigroups with a token $t(c_1,c_2)$ between
  them and $q(c_2,c_1)$ be its dual token. Let $f(c_1)$ be a
  combinatorial function with a corresponding recurrence operator $R$,
  which is defined by its kernel $r(c_1,c_1')$
  \begin{displaymath}
    [Rf](c_1)=\int_{C_1} r(c_1,c_1')\, f(c_1')\, dc_1'
  \end{displaymath}
  Then the generating function $\hat{f}(c_1)$ satisfies the
  equation
  \begin{equation}
    [\tilde{R}\hat{f}](c_2)= t(e_1,c_2),
  \end{equation}
  where $\tilde{R}$ is defined by the kernel
  \begin{equation}\label{eq:r-modif}
    \tilde{r}(c_2,c_2')=\int_{C_1} \int_{C_1}
    t(c_1'',c_2)\, r(c_1'',c_1')\, q (c_1',c_2')\, dc_1'\,dc_1''.
  \end{equation}
\end{thm}
\begin{proof}
  We start from an observation that
  \begin{equation}\label{eq:r-inter}
    \int_{C_2} \tilde{r}(c_2,c_2')\, t(c_1,c_2')\, dc_2'=
    \int_{C_1} t(c_1'',c_2)\, r(c_1'',c_1)\, dc_1'',
  \end{equation}
  which follows from application to both sides of~\eqref{eq:r-modif} the
  integral operator with kernel $t(c_1,c_2')$ and the identity
  \begin{displaymath}
    \int_{C_2} q(c_1',c_2')\, t(c_2',c_1)\, dc_2'=\delta(c_1',c_1).
  \end{displaymath}
  Then
  \begin{eqnarray}
    [\tilde{R}\hat{f}](c_2)&=&\int_{C_2} 
    \tilde{r}(c_2,c_2')\, \hat{f}(c_2')\, dc_2 \nonumber\\
    &=&\int_{C_2} \hat{r}(c_2,c_2') \int_{C_1} f(c_1)\, t(c_1,c_2') 
    \,dc_1\, dc_2 \nonumber\\
    &=&\int_{C_1} \int_{C_2} \hat{r}(c_2,c_2')\, t(c_1,c_2')\, 
    dc_2\, f(c_1) \,dc_1 \nonumber\\
    &=&\int_{C_1} \int_{C_1} t(c_1'',c_2)\, r(c_1'',c_1)\, dc_1''\, 
    f(c_1) \,dc_1 \label{eq:r-use}\\
    &=&\int_{C_1} t(c_1'',c_2)  \int_{C_1} r(c_1'',c_1)\, f(c_1) \,dc_1\, 
    dc_1''
    \nonumber\\
    &=&\int_{C_1} t(c_1'',c_2)  \delta (c_1'')\, dc_1'' \nonumber\\
    &=& t(e_1,c_2) \nonumber
  \end{eqnarray}
  where we deduce~\eqref{eq:r-use} from~\eqref{eq:r-inter}.
\end{proof}
In fact, our calculation just says that $\tilde{R}g=T\delta$ for
$g=Tf$ and operator $\tilde{R}=TRT^{-1}$ if $Rf=\delta$. A touch of
non-triviality appears when we use special properties of functions
$t(c_1,c_2)$ and $q(c_2,c_1)$ be a pair of dual tokens: 
\begin{cor}
  If under assumptions of Theorem~\ref{th:wisdom-teeth} the operator
  $R$ is shift invariant on $C_1$, i.e. has a kernel of the form
  $r(c_1,c_1')=r(c_1'^{-1}c_1)$, then the operator $\tilde{R}$ is also
  shift invariant on $C_2$ with a kernel of the form
  $\hat{r}(c_2,c_2')=r(c_2'^{-1}c_2)$.
\end{cor}
\begin{proof}
  This follows from the property of t-transform to map a shift
  invariant operator on $C_1$ to a shift invariant operator on $C_2$.   
\end{proof}
\begin{rem}
 The above corollary
  clearly indicates what type of generating functions is reasonable: 
  \emph{a simple generating function can be constructed with respect 
    to a token $t(c_1,c_2)$ such that the recurrence operator $R$ is
    shift-invariant with respect to it.}
\end{rem} 
As it usually occurs a simple fact may have interesting
realizations.
\begin{example}\label{ex:Fibonacci}
  Let $f_n$ be the Fibonacci numbers. Their known recursion
  $f_n=f_{n-1}+f_{n-2}$ should be stated in a more accurate way with
  the help of agreement~\eqref{eq:agreement} as follows:
  \begin{eqnordisp}
    f_n-f_{n-1}-f_{n-2}=\delta_{n,0}, \qquad n\geq 0.
  \end{eqnordisp}
  Then the recurrence operator is
  given by $R=I-S-S^2$, where $S$ is a backward shift $Sf(n)=f(n-1)$ on
  $\Space[+]{N}{}$, particularly $R$ \emph{is shift
    invariant} with respect to token $t(x,n)=x^n$. 
  The corresponding kernel is
  \begin{eqnordisp}
    r(n,i)=\delta_{n,i}-\delta_{n-1,i}-\delta_{n-2,i},
  \end{eqnordisp} 
  which is a function of $n-i$.  For the pair of dual
  tokens $t(x,n)=x^n$ and $q(i,y)=\delta^{(i)}(y)/i!$ we obtain a
  transformed kernel:
  \begin{eqnarray}
    \tilde{r}(x,y) &= & \sum_{n=0}^\infty \sum_{i=0}^\infty 
    t(x,n)\, r(n,i)\, q(i,y) \nonumber\\
    &= & \sum_{n=0}^\infty \sum_{i=0}^\infty 
    x^n \left(\delta_{n,i}-\delta_{n-1,i}
      -\delta_{n-2,i}\right) \frac{\delta^{(i)}(y)}{i!}\nonumber\\
    &= & \sum_{n=0}^\infty x^n (1-x-x^2) \frac{\delta^{(n)}(y)}{n!}
    \nonumber\\
    &= & (1-x-x^2)\sum_{n=0}^\infty x^n \frac{\delta^{(n)}(y)}{n!}  
    \label{eq:Fibonacci-rec}
  \end{eqnarray} 
  From the Taylor expansion the sum
  in~\eqref{eq:Fibonacci-rec} obviously represents an integral kernel
  of the identity operator, therefore the transformation $\tilde{R}$
  is the operator of multiplication by $1-x-x^2$ and the generating
  function $\hat{f}(x)$ should satisfy equation
  \begin{displaymath}
    \tilde{R}f(x)=(1-x-x^2)\hat{f}(x)=x^0=1.
  \end{displaymath}
  Thus 
  \begin{eqnordisp}
    \hat{f}(x)=\sum_{n=0}^\infty f_n\, x^n=\frac{1}{1-x-x^2}   
  \end{eqnordisp}
  Note that operator $\tilde{R}$ of
  multiplication by $(1-x-x^2)$ is ``\emph{shift
    invariant}''---it commutes with the operator of multiplication by
  $x$ (both these operators are certain convolutions under Fourier
  transform).

  An attempt to construct a transformation $\tilde{R}$ with respect to
  a similar pair of dual tokens $t_1(x,n)=x^n/n!$ and
  $q_1(i,y)=\delta^{(i)}(y)$ will not enjoy the above simplicity.
\end{example}
\begin{example}\label{ex:bell1}
  Let $B_n$ be the Bell numbers~\cite{Rota64a}
  or~\cite[\S~3]{RotaTay93a}. The known recursion for
  them~\cite[\S~3]{RotaTay93a} can again be restated as
  \begin{eqnordisp}[eq:bell-recurr]
    B_n - \sum_{k=0}^{n-1} \bin{n-1}{k} B_k=\delta_{n,0}, \qquad n \geq
    0.   
  \end{eqnordisp} 
  Such a recurrence identity is \emph{shift invariant}
  with respect to $W$ defined by the identity $Wf(n)=f(n-1)/n$ on
  $\Space[+]{N}{}$. In other words it is \emph{shift
    invariant} with respect to token $t(x,n)=x^n/n!$.
  The kernel $r(n,i)$ corresponding to recurrence operator $R$
  in~\eqref{eq:bell-recurr} is
  \begin{eqnordisp}
    r(n,i)=\delta_{n,i}- \sum_{k=0}^{n-1} \bin{n}{k} \delta_{k,i}.
  \end{eqnordisp} 
  We apply a transformation with respect to a pair of
  dual tokens $t(x,n)=x^n/n!$ and $q(i,y)=\delta^{(i)}(y)$:
  \begin{eqnarray}
    \tilde{r}(x,y) &= & \sum_{n=0}^\infty \sum_{i=0}^\infty 
    t(x,n)\, r(n,i)\, q(i,y) \nonumber\\
    &= & \sum_{n=0}^\infty \sum_{i=0}^\infty 
    \frac{x^n}{n!} \left( \delta_{n,i}- \sum_{k=0}^{n-1} 
      \bin{n-1}{k} \delta_{k,i} 
    \right)\delta^{(i)}(y)\nonumber\\
    &= & \sum_{n=0}^\infty \sum_{i=0}^\infty 
    \frac{x^n}{n!}\, \delta_{n,i}\, \delta^{(i)}(y)
      - \sum_{n=0}^\infty \sum_{i=0}^\infty\, 
    \frac{x^n}{n!}\, \sum_{k=0}^{n-1} 
      \bin{n-1}{k}\, \delta_{k,i}\, 
    \delta^{(i)}(y)\nonumber\\
    &= & \sum_{n=0}^\infty  
    \frac{x^n}{n!}\, \delta^{(n)}(y)
      - \sum_{n=0}^\infty  \sum_{i=0}^{n-1}
    \frac{x^n}{n!}\, \frac{(n-1)!}{i!\, (n-1-i)!}\,
    \delta^{(i)}(y)\nonumber\\
    &= & \sum_{n=0}^\infty  
    x^n\,\frac{\delta^{(n)}(y)}{n!} 
      - \sum_{n=0}^\infty  \sum_{i=0}^{n-1}
    \frac{1}{n}\, x^n\, \frac{\delta^{(i)}(y)}{i!\, (n-1-i)!}
    \nonumber\\
    &= & \sum_{n=0}^\infty x^n\,\frac{\delta^{(n)}(y)}{n!} 
      - \sum_{n=0}^\infty \frac{x}{n}  \sum_{i=0}^{n-1}
    \frac{x^{n-i-1}}{(n-1-i)!}\, \frac{x^i\,\delta^{(i)}(y)}{i!}
    \label{eq:bell-trans1} 
  \end{eqnarray} 
  One verifies that an integration of a function $f(y)$
  with the sum 
  \begin{eqnordisp}
    \sum_{i=0}^{n-1}
    \frac{x^{n-i-1}}{(n-1-i)!} \frac{x^i\, \delta^{(i)}(y)}{i!}   
  \end{eqnordisp}
  in \eqref{eq:bell-trans1} produces the $(n-1)$-th term in the Taylor
  expansion of the product $e^x f(x)$. Thus the entire expression
  in~\eqref{eq:bell-trans1} is a kernel of the operator:
  \begin{displaymath} 
    \tilde{R}=I- {\textstyle\int}e^x, 
  \end{displaymath}
  where $I$ is the identity operator and ${\textstyle\int}$ 
  is the operator of the anti-derivation
  fixed by the condition $[{\textstyle\int}g](0)=0$, particularly
  ${\textstyle\int} x^n/n!= x^{n+1}/(n+1)!$. So generating
  function $\hat{b}(x)$ should satisfy the equation
  $(I-{\textstyle\int}e^x)\hat{b}(x)=x^0=1$.  Taking the derivative
  from both sides we found the differential equation
  $(D-e^x)\hat{b}(x)=0$. With the obvious initial condition
  $\hat{b}(0)=1$ it determines that $\hat{b}(x)=\exp(\exp(x)-1)$.
\end{example}

One may interpret~\eqref{eq:recurrent} as if $f(c_1)$ is a fundamental
solution to operator $R$. This is indeed a right way of thinking if
one can found a group action such that $R$ is a shift invariant. Then
the convolution $h(c_1)=[g*f](c_1)$ will give a solution to the
equation $[Rh](c_1)=g(c_1)$.

\section{Models for the Umbral Calculus}\label{sec:models-umbr-calc}
We present realizations (models) for different descriptions of the
umbral calculus mentioned in Section~\ref{se:umbral}.

\subsection{Finite Operator Description: a Realization}
As was mentioned in Subsection~\ref{ss:operator-d} the main
ingredients of the approach are polynomial sequence of binomial type
and shift invariant operators. As was already pointed in
Examples~\ref{ex:polin1}, \ref{ex:polin2} and \ref{ex:polin3} the
notion of token is a useful refinement of polynomial sequence of
binomial type. It is much more useful if we consider it together with
shift invariant operators on c-semigroups not just groups.  Next
Subsections contains some details of this.

\subsection{Hopf Algebra Description: a Realization}
We describe some applications of c-semigroups and tokens. First we can
add new lines to Proposition~\ref{pr:isomorphism}:
\begin{prop}
  Let $C_2$ will be a c-semigroup with a right source $e$. Let also
  $C_1$ be another c-semigroup and $t(c_1,c_2)$ be a token between them.
  Then the list in the Proposition~\ref{pr:isomorphism} for $C_2$
  can be extended via t-transform by the
\begin{enumerate}
\item[\textup{4.}] The subspace of convolutions over $C_1$.
\item[\textup{5.}] The subspace of linear functional over
  $\FSpace{C}{}(C_1)$.
\end{enumerate}
\end{prop}
\begin{proof}
  By Theorem~\ref{th:t-transform} the space of convolutions over $C_2$
  can be mapped to the space of convolutions over $C_1$ via
  t-transform. A backward mapping from convolutions over $C_1$ to
  convolutions over $C_2$ is given by the dual token $q(c_2,c_1)$ of
  associated delta family $Q(c_1)$.
\end{proof}
We establish the correspondence of four mentioned sets as linear
spaces. However two of them (convolutions over $C_1$ and $C_2$) have
the isomorphic structures as convolution algebras. This allows us to
transfer also their multiplication law on the space of linear
functionals.
\begin{defn}
  The product $l=l_1*l_2$ of two functionals $l_1$ and $l_2$ is again a
  linear functional corresponding to the convolution $S=S_1*S_2$,
  where convolutions $S_1$ and $S_2$ correspond to $l_1$ and $l_2$. It 
  acts on a function $p(c_1)$ as follows:
  \begin{displaymath}
    \scalar{l}{p(c'_1,\cdot)}=\int_{C_1} \scalar{l_1}{p(c_1,\cdot)} 
    \scalar{l_2}{p([c_1^{-1}c_1'],\cdot)}\,dc_1
  \end{displaymath}
\end{defn}
\begin{example}\label{ex:polin3}
  We return to the notations of Examples~\ref{ex:polin1}
  and~\ref{ex:polin2}. For a linear functional $l$ on
  $C_2=\Space{R}{}$ we introduce its t-transform $[\oper{T}l](n)$ as its
  action (with respect to $x$ variable) on the token $t(n,x)=p_n(x)$
  given by a sequence of polynomial type $p_n(x)$:
  \begin{displaymath}
    [\oper{T}l](n)=lT(n,x)=lp_n(x).
  \end{displaymath}
  Then the product $l$ of two functionals $l_1$ and $l_2$ are defined by
  the convolution on $\Space[+]{N}{}$:
  \begin{equation}\label{eq:x-hopf}
    [\oper{T}l](n)=\int_{\Space[+]{N}{}} l_1(k)\, l_2([k^{-1}n])\, dn =
    \sum_{k=0}^n l_1(k)\, l_2(n-k).
  \end{equation}
\end{example} 
The above definition of product can be found
in~\cite{Rota78}. The c-semigroup algebra over $\Space[+]{N}{}$ can
be naturally realized as the multiplicative algebra of a formal power
series of one variable $t$ (see Proposition~\ref{pr:series}).  So we
obtain one more face an umbral algebra as follows:
\begin{cor} \textup{\cite{Rota78}}
  The algebra of linear functionals on \Space{R}{} with the Hopf
  multiplication~\eqref{eq:x-hopf} is in a natural correspondence with
  the algebra of formal power series in one variable $t$
  \begin{displaymath}
    \{l_n\}_{n=0}^{\infty} \mapsto \sum_{n=o}^{\infty} l_n t^{n}.
  \end{displaymath}
\end{cor}

\subsection{The Semantic Description: a Realization} 

We start from the c-semigroup $S$, which is the direct product of many
copies of the c-semigroup $\Space[+]{N}{}$. The different copies of
$\Space[+]{N}{}$ are labelled by letters of an alphabet $A$, which is
denoted by Greek letters. An element $k$ of a copy $\Space[+]{N}{}$
labelled by a letter $\alpha$ can be written interchangeably as
$k_\alpha$ or $\alpha^k$. In this sense any function $p(k_\alpha,
l_\beta, \ldots)$ in $\FSpace{C}{0}(S)$ can be identified with the
polynomial
\begin{displaymath}
  \textsf{p}(\alpha, \beta,\ldots )=\sum_{k_\alpha, l_\beta, \dots }
  p(k_\alpha, l_\beta, \ldots) \alpha^k \beta^l \cdots.
\end{displaymath}
More over this identification send convolution $p*q$ of two functions
$p, q\in \FSpace{C}{0}(S)$ over $S$ to the product of polynomials
$\textsf{p}\textsf{q}$. So we will not distinguish a function
$p(k_\alpha, l_\beta, \ldots)$ and the corresponding polynomial
$\textsf{p}(\alpha, \beta,\ldots )$ any more as well as the
convolution algebra $\FSpace{C}{0}(S)$ and polynomial ring $D[A]$ over
alphabet $A$.

A linear functional $\object{eval}$ is defined as
\begin{displaymath}
  \object{eval}(p)=\int_S p(k_\alpha, l_\beta, \ldots)\, ds
\end{displaymath}
via a measure $ds$ on $S$ such that $ds(1_\alpha )=1$ for any $\alpha
\in A$. Because every $p$ has a compact (i.e., finite) support in $S$
the integration does not generate any difficulties. The
equality~\eqref{eq:eval-sep} express the fact that $S$ is the direct
sum of $\Space[+]{N}{} $'s and the measure $ds$ is the direct product
of corresponding measures $ds_\alpha $.

In spirit of \ref{it:augmentation} we require that the measure
$ds_\epsilon$ on the copy $\Space[+]{N}{}$ labelled by the
distinguished umbra $\epsilon$ is defined by the sequence
$(1,0,0,0,\ldots)$.

Now $\object{eval}$ defines an equivalence relation---umbral
equivalence on $\FSpace{C}{0}(S)$, namely $f\simeq g$,
$f,g\in\FSpace{C}{0}(S)$ if $(f-g)$ belongs to the kernel of
$\object{eval}$ (i.e., $\object{eval}(f-g)=0$). Because
$\object{eval}$ is additive but not multiplicative we see that umbral
equivalence is invariant under addition but multiplication (unlike the
support of functions belongs to different components in the direct
product of $\Space[+]{N}{}$).

Finally we will define $n.\beta$ as a function
\begin{displaymath}
  \beta^1\otimes \beta^0 \otimes \cdots \otimes \beta^0 +
  \beta^0\otimes \beta^1 \otimes \cdots \otimes \beta^0 +
  \beta^0\otimes \beta^0 \otimes \cdots \otimes \beta^1 
\end{displaymath}
in the $n$-th tensor power $\otimes_{k=1}^n \Space[+]{N}{}(\beta)$ of a
copy of $\Space[+]{N}{}$ labeled by $\beta$ (this means that
$\object{eval}$ has the identical distribution on all of them).
\comment{If we use comultiplication $\Delta:
  \FSpace{C}{}(\Space[+]{N}{})\rightarrow
  \FSpace{C}{}(\Space[+]{N}{})\otimes \FSpace{C}{}(\Space[+]{N}{})$
  defined as $\Delta \beta=\beta\otimes 1 +1\otimes \beta$ then we
  can see that $n.\beta$ is a \emph{synonymous} to $\Delta^n
  \beta$.}

This gives a realization for the classical umbral calculus described
via semantic approach in~\cite{RotaTay93a,RotaTay94a}.

\renewcommand{\thesection}{}
\section{Acknowledgments}

This paper was inspired by works of \person{Gian-Carlo Rota} and his
personal encouragement and interest supported me during its
preparation.

Prof.~\person{V.~Kharchenko} kindly enlightened me in the field of
algebra. I am grateful to anonymous referees who took a hard job to
read the initial version of the paper and suggested innumerable
improvements and corrections.  Last but not least the Editors of
\textsf{Zitschrift f\"ur Analysis und ihre Anwendungen} shown unusual
patience and persistence during reviewing process of this paper.

The essential part of the paper was written while  
author's stay at the Department of Mathematical Analysis, 
University of Ghent, as a Visiting Postdoctoral Fellow 
of the FWO-Vlaanderen (Fund of Scientific 
Research-Flanders), Scientific Research Network 3GP03196 
``Fundamental Methods and Technique in Mathematics'', 
whose hospitality I gratefully acknowledge.

\bibliographystyle{plain}
\bibliography{abbrevmr,akisil,analyse,acombin,arare}

\end{document}
